\documentstyle[12pt]{article}

\expandafter\ifx\csname amssym.def\endcsname\relax \else\endinput\fi
\expandafter\edef\csname amssym.def\endcsname{%
       \catcode`\noexpand\@=\the\catcode`\@\space}
\catcode`\@=11
\def\undefine#1{\let#1\undefined}
\def\newsymbol#1#2#3#4#5{\let\next@\relax
 \ifnum#2=\@ne\let\next@\msafam@\else
\ifnum#2=\tw@\let\next@\msbfam@\fi\fi \mathchardef#1="#3\next@#4#5}
\def\mathhexbox@#1#2#3{\relax
 \ifmmode\mathpalette{}{\m@th\mathchar"#1#2#3}%
\else\leavevmode\hbox{$\m@th\mathchar"#1#2#3$}\fi}
\def\hexnumber@#1{\ifcase#1 0\or 1\or 2\or 3\or 4\or 5\or 6\or 7\or 8\or
 9\or A\or B\or C\or D\or E\or F\fi}

\ifcase\@ptsize
  \font\tenmsa=msam10 \font\sevenmsa=msam7 \font\fivemsa=msam5
\or
  \font\tenmsa=msam10 scaled \magstephalf \font\sevenmsa=msam7 scaled
\magstephalf \font\fivemsa=msam5 scaled \magstephalf
\or
  \font\tenmsa=msam10 scaled \magstep1 \font\sevenmsa=msam7 scaled
\fi

\newfam\msafam
\textfont\msafam=\tenmsa
\scriptfont\msafam=\sevenmsa
\scriptscriptfont\msafam=\fivemsa
\edef\msafam@{\hexnumber@\msafam}
\mathchardef\dabar@"0\msafam@39
\def\dashrightarrow{\mathrel{\dabar@\dabar@\mathchar"0\msafam@4B}}
\def\dashleftarrow{\mathrel{\mathchar"0\msafam@4C\dabar@\dabar@}}

\def\ulcorner{\delimiter"4\msafam@70\msafam@70 }
\def\urcorner{\delimiter"5\msafam@71\msafam@71 }
\def\llcorner{\delimiter"4\msafam@78\msafam@78 }
\def\lrcorner{\delimiter"5\msafam@79\msafam@79 }
\def\yen{{\mathhexbox@\msafam@55 }}
\def\checkmark{{\mathhexbox@\msafam@58 }}
\def\circledR{{\mathhexbox@\msafam@72 }}
\def\maltese{{\mathhexbox@\msafam@7A }}

\ifcase\@ptsize
  \font\tenmsb=msbm10 \font\sevenmsb=msbm7 \font\fivemsb=msbm5
\or
  \font\tenmsb=msbm10 scaled \magstephalf \font\sevenmsb=msbm7 scaled
\magstephalf \font\fivemsb=msbm5 scaled \magstephalf
\or
  \font\tenmsb=msbm10 scaled \magstep1 \font\sevenmsb=msbm7 scaled
\fi

\newfam\msbfam
\textfont\msbfam=\tenmsb
\scriptfont\msbfam=\sevenmsb
\scriptscriptfont\msbfam=\fivemsb
\edef\msbfam@{\hexnumber@\msbfam}
\def\Bbb#1{{\fam\msbfam\relax#1}}
\def\widehat#1{\setbox\z@\hbox{$\m@th#1$}%
 \ifdim\wd\z@>\tw@ em\mathaccent"0\msbfam@5B{#1}%
\else\mathaccent"0362{#1}\fi}
\def\widetilde#1{\setbox\z@\hbox{$\m@th#1$}%
 \ifdim\wd\z@>\tw@ em\mathaccent"0\msbfam@5D{#1}%
\else\mathaccent"0365{#1}\fi}

\ifcase\@ptsize
  \font\teneufm=eufm10 \font\seveneufm=eufm7 \font\fiveeufm=eufm5
\or
  \font\teneufm=eufm10 scaled \magstephalf \font\seveneufm=eufm7
scaled \magstephalf \font\fiveeufm=eufm5 scaled \magstephalf
\or
  \font\teneufm=eufm10 scaled \magstep1 \font\seveneufm=eufm7 scaled
\fi

\newfam\eufmfam
\textfont\eufmfam=\teneufm
\scriptfont\eufmfam=\seveneufm
\scriptscriptfont\eufmfam=\fiveeufm
\def\frak#1{{\fam\eufmfam\relax#1}}

\csname amssym.def\endcsname

\expandafter\ifx\csname pre amssym.tex at\endcsname\relax \else \endinput\fi
\expandafter\chardef\csname pre amssym.tex at\endcsname=\the\catcode`\@
\catcode`\@=11
\newsymbol\boxdot 1200
\newsymbol\boxplus 1201
\newsymbol\boxtimes 1202
\newsymbol\square 1003
\newsymbol\blacksquare 1004
\newsymbol\centerdot 1205
\newsymbol\lozenge 1006
\newsymbol\blacklozenge 1007
\newsymbol\circlearrowright 1308
\newsymbol\circlearrowleft 1309
\undefine\rightleftharpoons
\newsymbol\rightleftharpoons 130A
\newsymbol\leftrightharpoons 130B
\newsymbol\boxminus 120C
\newsymbol\Vdash 130D
\newsymbol\Vvdash 130E
\newsymbol\vDash 130F
\newsymbol\twoheadrightarrow 1310
\newsymbol\twoheadleftarrow 1311
\newsymbol\leftleftarrows 1312
\newsymbol\rightrightarrows 1313
\newsymbol\upuparrows 1314
\newsymbol\downdownarrows 1315
\newsymbol\upharpoonright 1316
 
\newsymbol\downharpoonright 1317
\newsymbol\upharpoonleft 1318
\newsymbol\downharpoonleft 1319
\newsymbol\rightarrowtail 131A
\newsymbol\leftarrowtail 131B
\newsymbol\leftrightarrows 131C
\newsymbol\rightleftarrows 131D
\newsymbol\Lsh 131E
\newsymbol\Rsh 131F
\newsymbol\rightsquigarrow 1320
\newsymbol\leftrightsquigarrow 1321
\newsymbol\looparrowleft 1322
\newsymbol\looparrowright 1323
\newsymbol\circeq 1324
\newsymbol\succsim 1325
\newsymbol\gtrsim 1326
\newsymbol\gtrapprox 1327
\newsymbol\multimap 1328
\newsymbol\therefore 1329
\newsymbol\because 132A
\newsymbol\doteqdot 132B
 
\newsymbol\triangleq 132C
\newsymbol\precsim 132D
\newsymbol\lesssim 132E
\newsymbol\lessapprox 132F
\newsymbol\eqslantless 1330
\newsymbol\eqslantgtr 1331
\newsymbol\curlyeqprec 1332
\newsymbol\curlyeqsucc 1333
\newsymbol\preccurlyeq 1334
\newsymbol\leqq 1335
\newsymbol\leqslant 1336
\newsymbol\lessgtr 1337
\newsymbol\backprime 1038
\newsymbol\risingdotseq 133A
\newsymbol\fallingdotseq 133B
\newsymbol\succcurlyeq 133C
\newsymbol\geqq 133D
\newsymbol\geqslant 133E
\newsymbol\gtrless 133F
\newsymbol\sqsubset 1340
\newsymbol\sqsupset 1341
\newsymbol\vartriangleright 1342
\newsymbol\vartriangleleft 1343
\newsymbol\trianglerighteq 1344
\newsymbol\trianglelefteq 1345
\newsymbol\bigstar 1046
\newsymbol\between 1347
\newsymbol\blacktriangledown 1048
\newsymbol\blacktriangleright 1349
\newsymbol\blacktriangleleft 134A
\newsymbol\vartriangle 134D
\newsymbol\blacktriangle 104E
\newsymbol\triangledown 104F
\newsymbol\eqcirc 1350
\newsymbol\lesseqgtr 1351
\newsymbol\gtreqless 1352
\newsymbol\lesseqqgtr 1353
\newsymbol\gtreqqless 1354
\newsymbol\Rrightarrow 1356
\newsymbol\Lleftarrow 1357
\newsymbol\veebar 1259
\newsymbol\barwedge 125A
\newsymbol\doublebarwedge 125B
\undefine\angle
\newsymbol\angle 105C
\newsymbol\measuredangle 105D
\newsymbol\sphericalangle 105E
\newsymbol\varpropto 135F
\newsymbol\smallsmile 1360
\newsymbol\smallfrown 1361
\newsymbol\Subset 1362
\newsymbol\Supset 1363
\newsymbol\Cup 1264
 
\newsymbol\Cap 1265
 
\newsymbol\curlywedge 1266
\newsymbol\curlyvee 1267
\newsymbol\leftthreetimes 1268
\newsymbol\rightthreetimes 1269
\newsymbol\subseteqq 136A
\newsymbol\supseteqq 136B
\newsymbol\bumpeq 136C
\newsymbol\Bumpeq 136D
\newsymbol\lll 136E
 
\newsymbol\ggg 136F
 
\newsymbol\circledS 1073
\newsymbol\pitchfork 1374
\newsymbol\dotplus 1275
\newsymbol\backsim 1376
\newsymbol\backsimeq 1377
\newsymbol\complement 107B
\newsymbol\intercal 127C
\newsymbol\circledcirc 127D
\newsymbol\circledast 127E
\newsymbol\circleddash 127F
\newsymbol\lvertneqq 2300
\newsymbol\gvertneqq 2301
\newsymbol\nleq 2302
\newsymbol\ngeq 2303
\newsymbol\nless 2304
\newsymbol\ngtr 2305
\newsymbol\nprec 2306
\newsymbol\nsucc 2307
\newsymbol\lneqq 2308
\newsymbol\gneqq 2309
\newsymbol\nleqslant 230A
\newsymbol\ngeqslant 230B
\newsymbol\lneq 230C
\newsymbol\gneq 230D
\newsymbol\npreceq 230E
\newsymbol\nsucceq 230F
\newsymbol\precnsim 2310
\newsymbol\succnsim 2311
\newsymbol\lnsim 2312
\newsymbol\gnsim 2313
\newsymbol\nleqq 2314
\newsymbol\ngeqq 2315
\newsymbol\precneqq 2316
\newsymbol\succneqq 2317
\newsymbol\precnapprox 2318
\newsymbol\succnapprox 2319
\newsymbol\lnapprox 231A
\newsymbol\gnapprox 231B
\newsymbol\nsim 231C
\newsymbol\ncong 231D
\newsymbol\diagup 231E
\newsymbol\diagdown 231F
\newsymbol\varsubsetneq 2320
\newsymbol\varsupsetneq 2321
\newsymbol\nsubseteqq 2322
\newsymbol\nsupseteqq 2323
\newsymbol\subsetneqq 2324
\newsymbol\supsetneqq 2325
\newsymbol\varsubsetneqq 2326
\newsymbol\varsupsetneqq 2327
\newsymbol\subsetneq 2328
\newsymbol\supsetneq 2329
\newsymbol\nsubseteq 232A
\newsymbol\nsupseteq 232B
\newsymbol\nparallel 232C
\newsymbol\nmid 232D
\newsymbol\nshortmid 232E
\newsymbol\nshortparallel 232F
\newsymbol\nvdash 2330
\newsymbol\nVdash 2331
\newsymbol\nvDash 2332
\newsymbol\nVDash 2333
\newsymbol\ntrianglerighteq 2334
\newsymbol\ntrianglelefteq 2335
\newsymbol\ntriangleleft 2336
\newsymbol\ntriangleright 2337
\newsymbol\nleftarrow 2338
\newsymbol\nrightarrow 2339
\newsymbol\nLeftarrow 233A
\newsymbol\nRightarrow 233B
\newsymbol\nLeftrightarrow 233C
\newsymbol\nleftrightarrow 233D
\newsymbol\divideontimes 223E
\newsymbol\varnothing 203F
\newsymbol\nexists 2040
\newsymbol\Finv 2060
\newsymbol\Game 2061
\newsymbol\mho 2066
\newsymbol\eth 2067
\newsymbol\eqsim 2368
\newsymbol\beth 2069
\newsymbol\gimel 206A
\newsymbol\daleth 206B
\newsymbol\lessdot 236C
\newsymbol\gtrdot 236D
\newsymbol\ltimes 226E
\newsymbol\rtimes 226F
\newsymbol\shortmid 2370
\newsymbol\shortparallel 2371
\newsymbol\smallsetminus 2272
\newsymbol\thicksim 2373
\newsymbol\thickapprox 2374
\newsymbol\approxeq 2375
\newsymbol\succapprox 2376
\newsymbol\precapprox 2377
\newsymbol\curvearrowleft 2378
\newsymbol\curvearrowright 2379
\newsymbol\digamma 207A
\newsymbol\varkappa 207B
\newsymbol\Bbbk 207C
\newsymbol\hslash 207D
\undefine\hbar
\newsymbol\hbar 207E
\newsymbol\backepsilon 237F
\catcode`\@=\csname pre amssym.tex at\endcsname
\def\Box{\hbox{\vrule height1ex\kern-0.4pt
\vbox to 1ex{\hrule width1ex\vfil\hrule width1ex}\kern-0.4pt\vrule height1ex}}
\newcommand{\sqr}[2]{{{\vcenter{\vbox{\hrule height.#2pt
\hbox{\vrule width.#2pt height#1pt \kern#1pt
\vrule width.#2pt}
\hrule height.#2pt}}}}}

\newcommand{\ovl}{\overline}

\newcommand{\be}{\begin{equation}}

\newcommand{\ee}{\end{equation}}
\newcommand{\al}{\alpha}

\newcommand{\gm}{\gamma}

\newcommand{\kp}{\kappa}
\newcommand{\lm}{\lambda}

\newcommand{\rh}{\rho}
\newcommand{\sg}{\sigma}

\newcommand{\ph}{\phi}

\newcommand{\phv}{\varphi}
\newcommand{\ch}{\chi}
\newcommand{\Sg}{\Sigma}
\newcommand{\ps}{\psi}

\newcommand{\om}{\omega}

\newcommand{\raw}{\rightarrow}
\newcommand{\A}{\frak A}
\newcommand{\B}{\frak B}
\newcommand{\g}{\frak g}

\newcommand{\bib}{\bibitem}
\newcommand{\cin}{C^{\infty}}
\renewcommand{\S}{\mbox{$\cal S$}}

\renewcommand{\H}{\mbox{$\cal H$}}

\newcommand{\n}{\parallel}

\newcommand{\R}{{\Bbb R}}

\newcommand{\notp}{p \kern-.48em /}
\newcommand{\ci}{\cite}
\newcommand{\bea}{\begin{eqnarray}}
\newcommand{\eea}{\end{eqnarray}}
\newcommand{\ot}{\otimes}
\newcommand{\half}{\mbox{\footnotesize $\frac{1}{2}$}}

\topmargin = - 0.5 cm
\newcommand{\la}{\langle} 
\newcommand{\ra}{\rangle}
\newcommand{\No}{{\cal N}_0}
\renewcommand{\O}{{\cal O}}
\begin{document}
 \setlength{\baselineskip}{1.5\baselineskip}
\thispagestyle{empty}
\title{ The quantization of constrained systems: from symplectic
reduction to Rieffel induction\thanks{To appear in {\em Proc.\ of the
XIV'th Workshop on Geometric Methods in Physics, Bia\l owie\.{z}a,
1995}, eds. J.P. Antoine et.\ al., 1996}}
\author{ N.\ P.\ Landsman\thanks{ E.P.S.R.C.\ Advanced Research Fellow}\\
  Department of Applied Mathematics and Theoretical Physics\\
University of Cambridge\\ Silver Street, Cambridge CB3 9EW, U.K. \\
\mbox{} \hfill\\ DAMTP-96-05}
\date{\today}
\maketitle
 \begin{abstract} This is an introduction to the author's recent work
on constrained systems. Firstly, a generalization of the
Marsden-Weinstein reduction procedure in symplectic geometry is
presented - this is a reformulation of ideas of Mikami-Weinstein and
Xu. Secondly, it is shown how this procedure is quantized by Rieffel
induction, a technique in operator algebra theory.  The essential
point is that a symplectic space with generalized moment map is
quantized by a pre- (Hilbert) $C^*$-module. The connection with
Dirac's constrained quantization method is explained.  Three examples
with a single constraint are discussed in some detail: the reduced
space is either singular, or defined by a constraint with incomplete
flow, or unproblematic but still interesting. In all cases, our
quantization procedure may be carried out. Finally, we re-interpret
and generalize Mackey's quantization on homogeneous spaces. This
provides a double illustration of the connection between $C^*$-modules
and the moment map.  \end{abstract}
\newpage
\section{Generalized Marsden-Weinstein reduction}
Symplectic reduction is a technique to construct new symplectic
manifolds from old ones.  The subject is full of subtleties, and we
refer to \ci{AGJ} for a particularly careful discussion.  The main
idea is to start from a symplectic manifold $(\Sg,\om)$, select a
subspace $C$ (the constraint hypersurface, which is embedded in $\Sg$
by the identity map $i$), and quotient $C$ by the null foliation
${\cal F}_0$ defined by the null directions of the induced symplectic
form $i^*\om$ on $C$, to the extent that these are tangent to $C$. In
favourable circumstances the quotient $C/{\cal F}_0$ is a symplectic
manifold, referred to as the reduced space of $\Sg$ (with respect to
$C$).  A related construction is possible when $\Sg$ is merely a
Poisson manifold \ci{MR}, and the construction below generalizes to
that case; for simplicity, we restrict ourselves to the symplectic
case in what follows. In any case, classical reduction is a two-step
procedure: firstly the constraints are implemented (corresponding to
the choice of $C$); secondly, roughly speaking, gauge-equivalent
points are identified. As we shall see, in quantum theory only one of
these steps has to be taken. Dirac's approach to constrained
quantization \ci{Dir} selects the first step, whereas the author's
singles out the second.

Marsden-Weinstein reduction (cf.\ \ci{AM} and refs.\ therein) is a
special case of the above construction. Let $(S,\om_S)$ be a
symplectic manifold on which a connected Lie group $G$ (whose Lie
algebra we denote by $\g$, with dual $\g^*$) acts from the right in
strongly Hamiltonian fashion. If $(\g^*)^-$ stands for $\g^*$ equipped
with {\em minus} the Lie-Poisson structure \ci{AM,MR}, one obtains an
equivariant moment map $J:S\raw (\g^*)^-$ which is a Poisson morphism.
Putting $J_X=\la J,X\ra$, this means that $\{J_X,J_Y\}=-J_{[X,Y]}$,
which is equivalent to the global property $J(xg)={\rm Ad}^*(g)J(x)$.
Choose a coadjoint orbit $\O\in\g^*$, and put $\Sg=S\times {\O}$, with
symplectic form $\om=\om_S+\om_{\O}$, where $\om_{\O}$ is the Kirillov
(etc.) symplectic form on the orbit (note that this form is the one
induced on $\O$ regarded as a symplectic leaf of $\g^*$ equipped with
{\em plus} the Lie-Poisson structure). The inclusion map of $\O$ into
$\g^*$ is called $\rh$. The constraint hypersurface is given by the
fiber product $C=S*_{\g^*}\O$, which stands for the collection of
points $\{(x,y)\in S\times \O|J(x)=\rh(y)\}$ (which is evidently
diffeomorphic to $J^{-1}(\O)$). If $\O$ consists of regular values of
the moment map, $C$ is co-isotropically immersed in $\Sg$, and the
reduced space $S^{\O}=C/{\cal F}_0$ coincides with the
Marsden-Weinstein quotient $C/G\simeq J^{-1}(\O)/G$ (cf.\
\ci{AGJ,SL} for the singular case).

We now generalize Marsden-Weinstein reduction as follows \ci{Lan1}.
Instead of $\g^*$ we consider an arbitrary Poisson manifold $P$, with
realization $S$, i.e., we suppose a Poisson morphism $J:S\raw P^-$ is
given (the `generalized moment map'), with $S$ symplectic.  Let a
second realization $\rh:S_{\rh}\raw P$ be given. We may then repeat
the above steps: we can form the fiber product $C=S*_P S_{\rh}$,
which, under the assumption that either $J_*$ or $\rh_*$ is surjective
at all points of $S$ or $S_{\rh}$ that are relevant to $C$, turns out
to be co-isotropically immersed in $\Sg=S\times S_{\rh}$. The reduced
space is then defined by $$S^{\rh}=(S*_P S_{\rh})/{\cal F}_0,$$ as
before (the null foliation is generated by the Hamiltonian vector
fields of functions of the type $J^*f-\rh^*f$, $f\in C^{\infty}(P)$).
This is essentially the reduction procedure of Mikami-Weinstein
\ci{MiW} and Xu \ci{Xu91}, but reformulated without reference to
symplectic groupoids (also cf.\ \ci{Zak}).

Let ${\cal A}$ be the Poisson subalgebra of $C^{\infty}(S)$ consisting
of all functions $f$ with the property that for all $g\in
C^{\infty}(P)$ the Poisson bracket $\{J^*g,f\}$ vanishes at all points
of $S$ which are relevant to $C$. We can define an `induced
representation' $\pi^{\rh}$ of $\cal A$ on $S^{\rh}$ by
$\pi^{\rh}(f)([x,y])=f(x)$, where $[x,y]\in S^{\rh}$ is the image of
$(x,y)\in C$ under the canonical projection from $C$ to $C/{\cal
F}_0$. By our definition of $\cal A$ this is well-defined (i.e.,
independent of the choice of $x$ in $[x,y]$).  In the context of
constrained systems, such functions $f$ are called weak observables;
an observable is an equivalence class of such functions under the
equivalence relation $f_1\sim f_2$ iff
$\pi^{\rh}(f_1)=\pi^{\rh}(f_2)$.

Our thinking of symplectic reduction as an induction construction is
further supported by a classical imprimitivity theorem \ci{Xu91,Lan1}
and by a theorem on symplectic reduction in stages
\ci{Lan1,MR,SL}, both of which parallel the (much older!)
corresponding 
theorems on induced
representations in a Hilbert space and operator algebra context
\ci{Rie}.
 \section{Quantized Marsden-Weinstein reduction} The problem of
quantizing constrained (i.e., reduced) systems was first addressed by
Dirac \ci{Dir}.  He realized that only one of the two steps of the
classical reduction procedure needs to be implemented at the quantum
level. Let us illustrate his approach in the example of
Marsden-Weinstein reduction by a compact connected group at level zero
(that is, the classical reduced space is $S^0=J^{-1}(0)/G$ - we assume
that 0 is a regular value of the moment map so that the reduced space
is a manifold).  The classical phase space $S$ is supposed to be
quantized by some Hilbert space $\H$, the state space of the
unconstrained system. The strongly Hamiltonian action of $G$ on $S$ is
quantized by a unitary representation $U$ of $G$ on $\H$ (here and in
what follows, all unitary representations are assumed to be
continuous). This representation should be compatible with the
requirement that for each $X\in\g$, the function $J_X$ on $S$ is
quantized by $i dU(X)$. For reduction at zero, the classical
constraints are $J_{T_a}=0$, where $\{T_a\}_a$ is a basis of $\g$.
Dirac then proposed that the Hilbert space of the constrained quantum
system $\H^0_D$ be given by the subspace of $\H$ which satisfies the
quantum constraints, i.e., $$
\H^0_D=\{\ps\in\H|\, dU(T_a)\ps=0\: \forall a\}. 
$$ Hence (since $G$ is connected) $\H^0_D=\H_{\rm id}=P_{\rm id}\H$,
where $P_{\rm id}\H$ is the projector on the subspace of $\H$ which
transforms trivially under $U(G)$. A (weak) observable is a
self-adjoint operator $A$ (assumed to be bounded for simplicity) which
commutes with $P_{\rm id}$.  For such operators, the restriction
$P_{\rm id}AP_{\rm id}=AP_{\rm id}$ to $P_{\rm id}\H$ is well-defined,
and they evidently act on $\H^0_D$ by $\pi^0_D(A)=AP_{\rm id}$.

In the context of a specific quantization scheme, such as geometric
quantization, one may ask whether $\H^0_D$ (in conjunction with an
action of the observables on it) coincides with the quantization of
$S^0$ (in other words, do reduction and quantization commute). The
mathematical literature on this difficult question goes back to (at
least) \ci{GS}, and we will not concern ourselves with this issue
here. A discussion of examples involving compact as well as
non-compact groups, using the quantization scheme of \ci{Lan1} and the
present paper, may be found in \ci{LanH}.

Rather, we propose to look at the above construction in a different
light. Instead of implementing the constraints, we can mimick the
second step of classical reduction, viz.\ quotienting by the null
foliation of the symplectic form. The way to do this consists in
modifying the inner product on $\H$. We define a new sesquilinear form
$(\, ,\,)_0$ on $\H$ by $$ (\ps,\phv)_0=\int_G dg\,(U(g)\ps,\phv), $$
where $dg$ is the normalized Haar measure on $G$.  From standard
compact group theory, we infer that $(\ps,\phv)_0=(P_{\rm
id}\ps,P_{\rm id}\phv)$.  Hence the new form is positive
semi-definite, with null space $\No=\H_{\rm id}^{\perp}$.  (Here and
in what follows, the null space of a positive semi-definite
sesquilinear form $(\, ,\, )_0$ is defined as the collection of
vectors $\ps$ such that $(\ps,\ps)_0=0$; by the Cauchy-Schwartz
inequality - which is the same as for positive definite sesquilinear
forms - it follows that $(\ps,\phv)_0=0$ for all $\phv$.)  We now
define the induced space $\H^0$ as $\H/\No$; clearly, $\H^0\simeq
\H_{\rm id}=\H_D^0$ as vector spaces. Let $V\ps$ be the image of
$\ps\in\H$ in $\H/\No$; the latter becomes a Hilbert space under the
inner product $(V\ps,V\phv)^0=(\ps,\phv)_0$. Since we have precisely
eliminated the null space of $(\, ,\,)_0$, this inner product is
positive definite, and $\H^0\simeq \H^0_D$ also as Hilbert spaces.

 The condition for an operator $A$ on $\H$ to have a well-defined
quotient action on $\H^0$ is $A\No\subseteq \No$.  We can then define
the (trivially) induced representation $\pi^0$ by
$\pi^0(A)V\ps=VA\ps$.  The requirement $\pi^0(A)^*=\pi^0(A^*)$ is
equivalent to $(A\ps,\phv)_0=(\ps,A^*\phv)_0$, which in turn is the
condition $[A,P_{\rm id}]=0$ we already encountered.  The property
$(A\ps,\phv)_0=(\ps,A^*\phv)_0$ actually implies that $A\No\subseteq
\No$, so that we can define a (weak) observable as a self-adjoint
operator $A$ which satisfies $(A\ps,\phv)_0=(\ps,A\phv)_0$. Thus the
present reformulation is completely equivalent to the Dirac approach.

Now consider the case where the reduced space $S^{\O}$ is obtained by
Marsden-Weinstein reduction from a nontrivial coadjoint orbit $\O$.
Apart from $\H$, the quantization procedure then needs a second
Hilbert space $\H_{\rh}$, seen as the quantization of $\O$, and an
irreducible unitary representation $U_{\rh}$ of $G$, which quantizes
the coadjoint action of $G$ on $\O$. The quantization of $S\times \O$
equipped with given $G$-action on $S$ times the coadjoint action is
then given by $\H\ot\H_{\rh}$ equipped with the unitary representation
$U\ot U_{\rh}$.  However, since $S^{\O}\simeq (S\times \O)^0$ as
symplectic spaces (where the moment map on $S\times\O$ is the sum of
$J$ and the inclusion map), its quantization according to Dirac can be
inferred from the preceding as being $$\H^{\rh}_D= P^{\rh}_{\rm id}
(\H\ot\H_{\rh}),$$ where $P^{\rh}_{\rm id}$ projects onto the subspace
of $\H\ot\H_{\rh}$ transforming trivially under $U\ot U_{\rh}$.

For the purpose of widening the discussion to general constrained
systems, it is convenient to reformulate the construction of Dirac's
Hilbert space $\H^{\rh}_D$ as follows. The group algebra $C^*(G)$ is
defined as the $C^*$-completion of the convolution algebra $L^1(G)$,
with adjoint defined by $f^*(g)=\overline{f(g^{-1})}$ \ci{Dix}.  This
algebra can be `diagonalized' by the Plancherel transform \ci{Dix}:
for $\gm\in\hat{G}$ (the unitary dual of $G$, i.e., the space of
equivalence classes of irreducible unitary representations of $G$; for
$G$ compact this is a discrete space) we put $\hat{f}(\gm)=\int_G dg\,
f(g)U_{\gm}(g)$ (initially defined for $f\in L^1(G)$ and extended to
$C^*(G)$ by continuity), with inverse $f(g)=\sum_{\gm\in\hat{G}}
d_{\gm}{\rm Tr}\, [U_{\gm}(g)^*\hat{f}(\gm)] $ (where $d_{\gm}$ is the
dimension of $\H_{\gm}$). Note that $f(\gm)$ is an operator on
$\H_{\gm}$, and that some choice $U_{\gm}$ has been made for the
unitary representatives in the equivalence class $\gm$.  The algebra
$C^*(G)$ then consists of those operator-valued functions $\hat{f}$ on
$\hat{G}$ for which the function $\gm\raw \n \hat{f}(\gm)\n$ is in
$C_0(\hat{G})$ (i.e., vanishes at infinity). One has
$\widehat{f^*}(\gm)=\hat{f}(\gm)^*$ and
$\widehat{f_1*f_2}(\gm)=\hat{f_1}(\gm)\hat{f_2}(\gm)$.  The norm in
$C^*(G)$ is given by $\n f\n={\rm
sup}_{\gm\in\hat{G}}\n\hat{f}(\gm)\n$.

Given a unitary representation $U(G)$ on a Hilbert space $\H$, we
obtain a representation $\pi$ of $C^*(G)$ by $\pi(f)=\int_G dg\,
f(g)U(g)$.  We use the decomposition $\H\simeq
\oplus_{\ch\in\hat{U}}\H(\ch)$ under $U(G)$, where
$\hat{U}\subseteq\hat{G}$, and $\H(\ch)=\H_{\ch}\ot {\cal K}_{\ch}$,
the second factor taking account of the multiplicity of the
irreducible $U_{\ch}$ in $\H(\ch)$.  This yields $\pi(f)\simeq
\oplus_{\ch\in\hat{U}} \hat{f}(\ch)\ot{\Bbb I}_{{\cal K}_{\ch}}$.
(The symbol ${\Bbb I}_{{\cal K}}$ denotes the identity operator on a
Hilbert space $\cal K$.)  Similarly, we obtain a representation
$\pi_{\rh}(C^*(G))$ on $\H_{\rh}$.

 Also, a right-representation $\pi_R$ on $\H$ (that is,
$\pi_R(AB)=\pi_R(B)\pi_R(A)$) may be defined by $\pi_R(f)=\int_G dg\,
f(g)U(g^{-1})$. This may be decomposed as $\pi_R(f)\simeq
\oplus_{\ch\in\hat{U}} \hat{f}(\ovl{\ch})^T\ot{\Bbb I}_{{\cal
K}_{\ch}}$ By tensoring with the relevant identity operators, $\pi_R$
and $\pi_{\rh}$ are defined on $\H\ot\H_{\rh}$. It then follows that
$\H^{\rh}_D$ may alternatively be characterized as the subspace of
$\H\ot\H_{\rh}$ consisting of those $\Psi$ on which the `quantum
constraints' $$ (\pi_R(f) - \pi_{\rh}(f))\Psi=0 \:\:\:\: \forall f\in
C^*(G) $$ hold.  For we can decompose $\Psi\simeq
\oplus_{\ch\in\hat{U}}\Psi(\ch)$, where $ \Psi(\ch)\in
\H(\ch)\ot\H_{\rh}$. The quantum constraints imply that $$
(\hat{f}(\ovl{\ch})^T\ot{\Bbb I}_{{\cal K}_{\ch}}\ot{\Bbb
I}_{\H_{\rh}}-{\Bbb I}_{\H}
\ot{\Bbb I}_{{\cal K}_{\ch}}\hat{f}(\rh))\Psi(\ch)=0
\:\:\:\forall \ch\in\hat{U}.
$$ For $\ch\neq\ovl{\rh}$ we choose $f$ such that
$\hat{f}(\ovl{\ch})={\Bbb I}_{\H_{\ch}}$ and $\hat{f}(\rh)=0$, which
shows that $\Psi(\ch)=0$ for such $\ch$.  For $\ch = \ovl{\rh}$ the
constraints are only satisfied if $\Psi(\ovl{\rh})=
\sum_i e_i\ot\ps\ot e_i$, where $\{e_i\}$
 is a basis in $\H_{\rh}$ (and in $\H_{\ovl{\rh}}$), and
$\ps\in {\cal K}_{\ch}$ is arbitrary. But these vectors precisely form
$ P^{\rh}_{\rm id}\H\ot\H_{\rh}$.

The above way of writing $\H^{\rh}_D\subset \H\ot\H_{\rh}$ reflects
the fact that to some extent it is the quantum analogue of the
constraint hypersurface $C=S*_{\g^*}\O\subset S\times \O$ of the
classical theory.

Evidently, we may arrive at $\H^{\rh}_D$ by the second method of
modifying the inner product in putting $$ (\Psi,\Phi)_0=\int_G
dg\,(U(g)\ot U_{\rh}(g)\Psi,\Phi) = (P^{\rh}_{\rm id}\Psi,P^{\rh}_{\rm
id}\Phi).  $$ The null space is then given by $\No=(P^{\rh}_{\rm
id}\H\ot\H_{\rh})^{\perp}$, so that the induced space
$\H^{\rh}=\H\ot\H_{\rh}/\No$ (equipped with the inner product
inherited from $(\, ,\,)_0$) coincides with $\H^{\rh}_D$.

 Also in this case we may rewrite the construction, making use of
$C^*(G)$.  Consider the function $\la \ps,\phv\ra_{C^*(G)} $ on $G$
defined by $g\raw (U(g)\phv,\ps)$, where $\ps,\phv\in\H$.  This
defines a map $\la\cdot,\cdot\ra_{C^*(G)}:\overline{\H}\ot \H\raw
C^*(G)$, which behaves like an inner product taking values in the
$C^*$-algebra $C^*(G)$. Namely, one has
$\la\ps,\phv\ra_{C^*(G)}^*=\la\phv,\ps\ra_{C^*(G)}$. Moreover, it
enjoys the `equivariance' property $
\la\ps,\pi_R(f)\phv\ra_{C^*(G)}=\la\ps,\phv\ra_{C^*(G)}f$, where the
right-hand side contains the (convolution) product of two elements of
$C^*(G)$.  Finally, one has $\la\ps,\ps\ra_{C^*(G)}\geq 0$ for all
$\ps\in\H$.  To show this, expand $\ps=\sum_{\ch\in\hat{U}}
\ps_i(\ch)e_i(\ch)\ot\kappa(\ch)$, where $\{e_i\}$ is a basis in
$\H_{\ch}$, and $\kappa(\ch)\in {\cal K}_{\ch}$ has unit norm.  The
Plancherel transform of $\la\ps,\ps\ra_{C^*(G)}$ is then given by
$\gm\raw 0$ if $\gm\notin \hat{U}$, and $\gm\raw
\ps_i(\ovl{\gm})\ovl{\ps_j(\ovl{\gm})}$ if $\ovl{\gm}\in\hat{U}$.
Positivity in $C^*(G)$ simply means that $\hat{f}(\gm)$ is a positive
matrix for all $\gm$, which is clearly the case (cf.\ \ci{Lan1} for a
different proof).

The map $\la\cdot,\cdot\ra_{C^*(G)}$ is the quantum analogue of the
moment map
\ci{Lan1}.
  If $\ps\in\H(\ch)$ then the Plancherel transform of $\la\ps,\ps\ra_{C^*(G)}$
only has support at $\gm=\ovl{\ch}$, so that this `quantum moment map'
detects which representation of $G$ is carried by a vector $\ps$; this
reflects the fact in symplectic geometry that the $G$-orbit through a
point of $S$, in case this orbit happens to be symplectic, is carried
into a given co-adjoint orbit by the moment map.

We can now rewrite the modified inner product on $\H\ot\H_{\rh}$ by
linear extension of $$ (\ps\ot v,\phv\ot w)_0=(\pi_{\rh}(\la
\phv,\ps\ra_{C^*(G)})v,w)_{\rh}; $$ the inner product on the
right-hand side is the one in $\H_{\rh}$.  This is a sesquilinear form
by the `inner product' property shown 3 paragraphs ago.  It is
positive semi-definite by the positivity property proved afterwards.
Define $$D_0=\{(\pi_R(f)-\pi_{\rh}(f))\Psi|f\in
C^*(G),\Psi\in\H\ot\H_{\rh}\}.$$ A short computation shows that
$D_0\subseteq\No$.  Using the Plancherel transform in similar vein to
the argument leading to the alternative identification of
$\H^{\rh}_D$, one can show the opposite inclusion $\No\subseteq D_0$.
Hence $\No=D_0$. In other words, every vector satisfies the quantum
constraints up to vectors in the null space of $(\, ,\, )_0$.  Since
these vectors project to zero in the induced space $\H^{\rh}$, we
clearly see how the second method to quantize the constrained systems
in questions gets rid of the states not satisfying the constraints,
without actually having to impose the latter.

There are two reasons why the Dirac method was successful in this
class of examples: firstly, the constraints were first-class, and
secondly (because the group was compact) they had 0 in their discrete
spectrum (assuming that $\H_{\rm id}$ is not empty) with common
eigenspace. However, few realistic systems meet these conditions. If
there are second-class constraints (that is, constraints whose Poisson
bracket does not vanish on $C$), Dirac and his followers have to get
rid of them already at the classical level. If the first-class
constraints fail to have 0 as an eigenvalue, it is not clear what to
do (see below).  Fortunately, it turns out that the second approach
(of modifying the inner product) continues to work even in cases where
Dirac's method faces difficulties.
\section{Quantized symplectic reduction}
Returning to the general reduction procedure discussed at the end of
section 1, we now assume that, apart from the quantization of $S$ by
$\H$, we have quantized $S_{\rh}$ by a Hilbert space $\H_{\rh}$, and
the Poisson algebra $\cin(P)$ by a $\mbox{}^*$-algebra $\B$. For
technical reasons, we assume that $\B$ is a $C^*$-algebra or a
pre-$C^*$-algebra (this could correspond to the quantization of, say,
the bounded functions in $\cin(P)$); it may be thought of as the
abstract operator algebra generated by the quantum constraints.  The
pull-back $J^*:\cin(P)\raw\cin(S)$ is quantized by a
right-representation $\pi_R$ of $\B$ on $\H$. Also, the Poisson
morphism $\rh:S\raw P$ (or rather its pull-back) is quantized by a
representation $\pi_{\rh}$ of $\B$ on $\H_{\rh}$. As before, $\pi_R$
and $\pi_{\rh}$ may be regarded as operators on $\H\ot\H_{\rh}$ by
tensoring with the appropriate identity operators.

Recall that the classical constraint hypersurface $C$ was defined as
the fiber product $C=S*_P S_{\rh}\subset S\times S_{\rh}$, whose
points $(x,y)$ are singled out by the condition $J(x)=\rh(y)$.
Generalizing the expression for the quantum constraints in the compact
group case discussed above, Dirac's method would analogously attempt
to single out the subspace $\H^{\rh}_D$ of $\H\ot\H_{\rh}$ satisfying
$$
\pi_R(B)\Psi=\pi_{\rh}(B)\Psi\:\:\: \forall B\in \B.
$$ However, $\H^{\rh}_D$ is empty whenever at least one operator
$\pi_R(B)-\pi_{\rh}(B)$ fails to have zero in its discrete spectrum.
The empty space is not the correct quantization of the classical
reduced space $S^{\rh}$. (Even if $\H^{\rh}_D$ is nonempty, it may not
be the correct quantization.)  Hence we have to look for a different
approach.

Generalizing the construction of the previous section, our goal is to
define a sesquilinear form $(\, ,\,)_0$ on $\H\otimes\H_{\rh}$ with
certain properties. It turns out that very often the domain of
definition of this form must be taken to be a subspace $L\ot\H_{\rh}$,
where $L\subset\H$ is dense. In such cases, $(\, ,\,)_0$ may even be
non-closable as a quadratic form on $\H\otimes\H_{\rh}$, cf.\
\ci{LW1}. Also, if $L\neq \H$ we require that $\pi_R(B)$ maps $L$ into
itself. In fact, for the following construction is not even necessary
that $L$ is a subspace of a Hilbert space.
 
Firstly, $(\, ,\,)_0$ has to be positive semi-definite (i.e,
$(\Psi,\Psi)_0\geq 0$ for all $\Psi\in L\otimes\H_{\rh}$); as before,
we denote its null space by $\No$.  The `induced' space $\H^{\rh}$ is
now defined as the completion of $L\otimes\H_{\rh}/\No$ under the
(obvious) inner product defined by
$$(V\Psi,V\Phi)^{\rh}=(\Psi,\Phi)_0,$$ where $V:L\otimes\H_{\rh}\raw
L\otimes\H_{\rh}/\No$ is the canonical projection map. In other words,
the inner product on $\H^{\rh}$ is essentially $( \, ,\, )_0$ modulo
the null vectors.

 Secondly, all vectors of the form $(\pi_R(B) -\pi_{\rh}(B))\Psi$
should lie in $\No$.  As before, this means that every vector
satisfies the quantum constraints up to terms which project to the
zero vector in $\H^{\rh}$.

In a totally different context, an inner product meeting these
requirements was given by Rieffel
\ci{Rie}. Namely, if we can find a `generalized quantum moment map' 
$\la\cdot,\cdot\ra_{\B}:\overline{L}\ot L\raw \B$, with (i) the
`$\B$-valued inner product' property
$\la\ps,\phv\ra_{\B}^*=\la\phv,\ps\ra_{\B}$, (ii) the positivity
property $\pi_{\rh}(\la \ps,\ps\ra_{\B})\geq 0$ (as an operator on
$\H_{\rh}$) for all $\ps\in L$, and (iii) the `equivariance' property
$ \la\ps,\pi_R(B)\phv\ra_{\B}=\la\ps,\phv\ra_{\B}B$ for all $B\in\B$
and $\psi,\phv\in L$, then the form defined by linear extension of $$
(\ps\ot v,\phv\ot w)_0=(\pi_{\rh}(\la \phv,\ps\ra_{\B})v,w)_{\rh} $$
satisfies the two conditions. The three properties imply the
sesquilinearity, the first condition, and the second condition,
respectively.

In the mathematical literature (e.g., \ci{Con}) one finds the closely
related concept of a (Hilbert) $C^*$-module: this is given by the
above data $(L, \B, \pi_R, \la\cdot,\cdot\ra_{\B})$ if in addition
$\B$ is a $C^*$-algebra, $\la\ps,\ps\ra_{\B}\geq 0$ (as an element of
$\B$ - this of course implies the weaker positivity property we
imposed above), and $L$ is complete under the norm $\n \ps\n^2=\n
\la\ps,\ps\ra_{\B}\n$. In examples relevant to quantization theory one
usually does not have $C^*$-modules, not even when $L$ is a Hilbert
space. (For example, if $G=U(1)$ acts on $L=L^2(U(1))$ in the regular
representation, with our choice of the quantum moment map the norm on
$L$ as a $C^*$-module is $\n \ps\n ={\rm sup}_n |\ps_n|$, where the
$\ps_n$ are the Fourier coefficients of $\ps$. The completetion of $L$
in this norm is $C_0({\Bbb Z})$, which is strictly bigger than
$L\simeq L^2({\Bbb Z})$.)

Weak quantum observables $A$ of the constrained system are defined as
operators on $L$ which are self-adjoint with respect to $(\, , \,)_0$,
i.e., $(A\Psi,\Phi)_0=(\Psi,A\Phi)_0$ for all $\Psi,\Phi\in
L\ot\H_{\rh}$ (here $A$ is identified with $A\ot {\Bbb I}$).  For such
operators, the induced representation \ci{Rie} $\pi^{\rh}$, given by
$\pi^{\rh}(A)V\Psi=VA\Psi$, is well-defined. This is because $\No$ is
mapped into itself by weak observables. In addition, $\pi^{\rh}(A)$ is
symmetric on the domain $V(L\ot\H_{\rh})$ if $A$ is a weak observable.
A sufficient condition for $A$ to be a weak observable is that $\la
A\ph,\psi\ra_{\B}=\la
\ph,A\psi\ra_{\B}$ for all $\ps,\phv\in L$. 
 Observables are equivalence classes of weak observables, where
$A_1\sim A_2$ iff $\pi^{\rh}(A_1)=\pi^{\rh}(A_2)$.
 
Even if $L\subset \H$ and $A$ is a weak observable which is bounded as
an operator on $\H$, the operator $\pi^{\rh}(A)$ is not necessarily
bounded on $\H^{\rh}$; one needs the bound $(A\Psi,A\Psi)_0\leq C_A
(\Psi,\Psi)_0$ for all $\Psi\in L\ot\H_{\rh}$ and some constant $C_A$
to prove boundedness of $\pi^{\rh}(A)$.  The $C^*$-algebra of weak
observables is defined as the completion of the collection of weak
observables which are bounded in this sense, equipped with a norm
given by $\n A\n =\sqrt{C_A}$, for the smallest possible $C_A$ in the
above bound (to be precise, operators whose norm is zero are to be
quotiented away). A particularly favourable case occurs when $L=\H$,
and the collection of weak observables happens to be a $C^*$-algebra
under the operator norm. In that case, the semi-positivity of $(\, ,\,
)_0$ implies that $\pi^{\rh}$ of a weak observable is bounded
\ci{Rie}. Unfortunately, if $L$ is a proper subspace of $\H$ then the
set of bounded operators which map $L$ into itself is never complete
under the operator norm.

There is a slight reformulation of the construction of the induced
space $\H^{\rh}$.  Given the quantum moment map
$\la\cdot,\cdot\ra_{\B}$ and a state $\om$ on $\B$, we can define a
modified inner product directly on $L$ by
$(\ps,\phv)_0^{\om}=\om(\la\phv,\ps\ra_{\B})$. The induced space
$\H^{\rh}_{\om}$ is then defined as the completion of $L/\No^{\om}$,
where $\No^{\om}$ is the null space of $(\, ,\, )_0^{\om}$. A map
$V_{\om}:L\raw L/\No^{\om}$ and corresponding induced representation
$\pi^{\rh}_{\om}$ of the $C^*$-algebra of weak observables $\A$ on $L$
may then be defined as before.  An arbitrary unit vector
$v\in\H_{\rh}$ defines a vector state $\om_v$ on $\B$, given by
$\om_v(B)=(\pi_{\rh}(B)v,v)$.  For each such $\om_v$, the pair
$(\H^{\rh}_{\om_v},\pi^{\rh}_{\om_v}(\A))$ is unitarily equivalent to
$(\H^{\rh},\pi^{\rh}(\A))$.  (Note that, for any $C^*$-algebra $\A$,
one may take $L=\B=\A$ and $\la A,B\ra_{\B}=A^*B$; the induction
construction is then equivalent to the GNS construction \ci{Dix}.)

It is possible to regard $\H^{\rh}$ as a subspace of the algebraic
dual $L^*$ of $L$.  This embedding is not canonical, and depends on
the choice of a state $\om_v$ on $\B$.  Firstly, consider vectors in
$\H^{\rh}$ of the type $V\psi$. Such a vector defines a linear
functional on $L$ by $<V\psi,\phv>=(\ps,\phv)_0^{\om_v}$. This is
well-defined, since $(\ps,\phv)_0^{\om_v}=0$ if $\ps\in\No^{\om_v}$.
Subsequently, if $V\ps_n\raw \ch$ in $\H^{\rh}$ then one can define
$<\ch,\phv>=\lim_n (\ps_n,\phv)_0^{\om_v}$. This is well-defined,
since $V\ps_n\raw 0$ in $\H^{\rh}$ is equivalent to
$(\ps_n,\ps_n)_0^{\om_v}\raw 0$, which implies $<\ps_n,\phv>\raw 0$
for all $\phv\in L$ by the Cauchy-Schwartz inequality for positive
semi-definite forms.  Hence we may regard $V$ as a map from $L$ into
$L^*$. Vectors in $\No^{\om_v}$ are mapped into the zero functional.
This way of looking at the theory shows how the constrained
quantization formalism of \ci{Ash} (also cf.\ \ci{Mar}) is a special
case of our method. (Their method corresponds to taking
$\H_{\rh}={\Bbb C}$, and $L$ is taken to be a topological vector space
which is continuously embedded in $\H$, such that $L\subset \H\subset
L'$ forms a Gel'fand triplet; here $L'\subset L^*$ is the topological
dual.)
\section{The case of commuting constraints}
 If classically we have $n$ commuting constraints $J_i\in\cin(S)$,
such that the constraints generate a group action, the reduced phase
space is a Marsden-Weinstein quotient $S^0=J^{-1}(0)/G$ with respect
to the group $G=\Pi_{i=1}^n G_i$, where $G_i$ is $\R$ or $U(1)$. The
case of torus actions ($G_i=U(1)$ for all $i$) has been well-studied,
and many beautiful results are available (such as the convexity
theorem of Atiyah and Guillemin-Sternberg, cf.\ \ci{Gui}). The torus
action on $S$ is quantized by a unitary representation $U$ on a
Hilbert space $\H$, which decomposes as $\H\simeq
\oplus_{ l} \H_{ l}$, where $ l=(l_1,\ldots,l_n)\in \hat{G}=
{\Bbb Z}^n$, and $\H_{ l}$
carries the representation $U_{ l}(z_1,\ldots,z_n)=z_1^{l_1}\ldots
z_n^{l_n}$.  According to the theory above, the quantization of $S^0$
is $\H_0$; see \ci{GS,Gui} for more information on this case.

We now look at the opposite case $G_i=\R$ for all $i$. In the
corresponding quantum theory, the Hilbert space $\H$ will carry a
unitary representation $U(\R^n)$.  By the SNAG theorem (or the
complete von Neumann spectral theorem), one has the direct integral
decomposition $\H\simeq \int_{\hat{G}} d\mu(\lm)\,\H(\lm)$ (cf.\
\ci{Dix}), where $\hat{G}=\R^n$, and the Hilbert space $\H(\lm)$
carries the representation $a\raw \exp(i\lm a)$.  If $\mu(\lm=0)$ is
positive, and $\lm=0$ is not in the essential spectrum, then the
quantization of $S^0$ is $\H(0)$.  Assume, instead, that 0 lies in the
absolutely continuous part of $\mu$; below we assume that there is no
discrete spectrum).

 Let $W:\H\raw \int_{\hat{G}} d\mu(\lm)\,\H(\lm)$ be a unitary
transformation that diagonalizes $U$. We choose the dense subspace
$L\subset \H$ in such a way that it is contained in the space of
elements $\ps$ of $\H$ for which $W\ps$ (regarded as a cross-section
of the field $\{\H(\lm)\}$
\ci{Dix}) is continuous on $\hat{G}$. We now define the form $(\, ,\,)_0$ as in
the compact group case (but only on the domain $L$), and compute: $$
(\ps,\phv)_0=\int_G da\, (U(a)\ps,\phv)=(2\pi)^n
\left| \frac{d\mu}{d\lm}\right| (0)((W\ps)(0),(W\phv)(0))_{\H(0)}. $$
Hence we can define $V:\H\raw \H^0=\H(0)$ by $$ V\ps=\left((2\pi)^n
\left| \frac{d\mu}{d\lm}\right| (0) \right)^{1/2}(W\ps)(0);
$$ this satisfies $(V\ps,V\phv)=(\ps,\phv)_0$, and shows that the
induced space $\H^0$ indeed coincides with $\H(0)$.

We illustrate this procedure in three examples, each of which is of
special interest for a different reason. The examples correspond to
the classical constraint $$ H_{\kp}=\half(p_x^2+\kp e^{4x}-p_y^2
)=0,$$ defined on $S=T^*\R^2$ with canonical symplectic structure, and
described in canonical co-ordinates, where the parameter $\kp$ assumes
the values $0,1,-1$. This constraint comes from a certain
finite-dimensional approximation to the universe as described by the
general theory of relativity, cf.\ \ci{LanW} and refs.\ therein.

Firstly, take $\kp=0$. This case also emerges in the context of the
representation theory of the Poincar\'{e} group of a two-dimensional
space-time, cf.\ \ci{LW1}.  The constraint generates an action of
$\R$, viz.\ $(x,y,p_x,p_y)\raw (x+p_xt,y-p_yt,p_x,p_y)$, with moment
map $J=H_0$. However, $J_*$ fails to be surjective at all points of
the form $(x,y,0,0)$ (which we will refer to as `the singular
points'), at which $J_*$ vanishes. Hence $0$ is not a regular value of
$J$, and $J^{-1}(0)$ is not a submanifold of $S$.  The
Marsden-Weinstein reduced space $S^0=J^{-1}(0)/\R$ does not have a
constant dimension: if we look at $S^0$ as fibered over the subspace
$p_x=\pm p_y$ of $\R^2$, then the fiber above $(0,0)$ is
two-dimensional whereas at all other points it is one-dimensional.
This example is of interest partly because the general theory of
singular Marsden-Weinstein reduction \ci{SL}, based as it is on the
assumption that the group action is proper, does not directly apply
here (the $\R$-action is not proper precisely at the singular points).

 Nonetheless, $J^{-1}(0)$ is strongly co-isotropic and locally conical
in the sense of \ci{AGJ}, so that the Marsden-Weinstein quotient $S^0$
agrees with the `geometric' reduction defined by $J^{-1}(0)$.  As
shown in \ci{AGJ}, it is therefore possible to define a Poisson
algebra $\hat{C}^{\infty}(S^0)$ as the space of strong observables
equipped with the Poisson bracket inherited from $S$. Elements of
$\hat{C}^{\infty}(S^0)$ are functions $f\in\cin(S)$, restricted to
$J^{-1}(0)$, which satisfy $\{f,H\}=0$ on $J^{-1}(0)$. It follows that
such $f$ depends on $x$ and $y$ through the combination $xp_y+yp_x$. A
study of the hamiltonian flow on $S$ defined by such functions, and
therefore of the corresponding flow on $S^0$ obtained by projection,
shows that $S^0$ may be decomposed into five `symplectic leaves' (cf.\
\ci{MR} for this concept in the regular case, and \ci{SL} for the
singular case with proper group action): $p_x=p_y>0$, $p_x=p_y<0$,
$p_x=-p_y>0$, $p_x=-p_y<0$, and $p_x=p_y=0$. Any point in a given leaf
cannot leave the leaf under a Hamiltonian flow.

For proper group actions, it is shown in \ci{SL} that the symplectic
leaves of a singular quotient $J^{-1}(0)/G$ are the components of
$(J^{-1}(0)\cap S_{(H)})/G$, where $S_{(H)}$ is the stratum in $S$ of
orbit type $H\subseteq G$ (that is, the stability group of any point
in $S_{(H)}$ is conjugate to $H$). In our example, the singular points
form the stratum of orbit type $\R$, and the other four leaves
correspond to the components having orbit type $\{e\}$. Hence we have
the same situation as for proper group actions.

We now turn to quantization. A remarkable feature of our quantization
scheme is that it applies even if the classical reduced space is
singular. We quantize $S$ by $\H=L^2(\R^2)$, realized in position
space. The constraint is quantized by the closure of the operator $
\hat{H}_0=\half(-(\partial/\partial x)^2+ (\partial/\partial y)^2)$,
initially defined and essentially self-adjoint on
$C_c^{\infty}(\R^2)$, which of course has absolutely continuous
spectrum equal to $\R$. We choose $L$ to be the subspace (easily shown
to be dense) of $\H$ consisting of those functions $\ps$ whose Fourier
transform $\check{\ps}$ is in $\cin(\R^2)$ and satisfies $\check{\ps}(
0)=0$.  It is not necessary that $L$ contain the domain of the quantum
constraint (indeed, it does not), since what matters is the unitary
group generated by it, which is defined on all of $\H$.

  One finds $$ (\ps,\phv)_0=(2\pi)^{-1}\int\frac{dp}{2|p|}\,
[\check{\ps}(p,p)\ovl{\check{\phv}(p,p)}+\check{\ps}(p,-p)
\ovl{\check{\phv}(p,-p)}].
$$ We see that the induced space may be realized as $\H^0=\oplus_{\pm}
L^2(\R,dp/4\pi |p|)$; the map $V:L\raw
\H^0$ assumes the form 
$(V\ps)_{\pm}(p)=\check{\ps}(p,\pm p)$. Interestingly, we may write $$
(V\ps)_{\pm}(p)=(\ps,f_{\pm}(p;\cdot)), $$ where
$f_{\pm}(p;x,y)=\exp(-ip(x\pm y))$. For each $p$, $f_{\pm}(p;\cdot)$
is a (generalized) solution to the quantum constraint
$\hat{H}_0f_{\pm}(p;\cdot)=0$. These solutions do not lie in $\H$, yet
the expression $(\ps,f_{\pm}(p;\cdot))$ is well-defined for $\ps\in L$
(cf.\ the comment below).

There is a quantum analogue of four of the five strata of the
classical reduced space (the stratum of orbit type $\R$ is not
represented in the quantum theory).  Classical (weak) observables had
to be smooth functions of $p_x$, $p_y$, and $xp_y+yp_x$.  On $\H$
these are quantized in the Schr\"{o}dinger representation, and the
corresponding induced representatives on $\H^0$ are given by $p\oplus
p$, $p\oplus -p$, and $-i(pd/dp \oplus -pd/dp)$, respectively. Note
that all three are essentially self-adjoint on $\oplus_{\pm}[
C_c^{\infty}(\R^+) \oplus C_c^{\infty}(\R^-)]$. Now each copy of
$L^2(\R,dp/4\pi |p|)$ in $\H^0$ splits as a direct sum
$L^2(\R^+,dp/4\pi |p|)\oplus L^2(\R^-,dp/4\pi |p|)$, each summand of
which is irreducible under the unitary group generated by the Lie
algebra spanned by the operators in question; this construction
realizes the four inequivalent massless representations of the
two-dimensional Poincar\'{e} group.

A similar argument might be given in terms of the flow generated by
the quantum observables.  By definition, the latter are induced
representatives $\pi^0(A)$ of weak observables $A$ (which, in
particular, must map $L$ into itself).  The difficulty is that
realistic Hamiltonians are unbounded operators, so that one has to
answer the question whether essential self-adjointness of $A$ on $L$
implies the same for $\pi^0(A)$ on $V(L)$, and if so, whether the
restriction that $\pi^0(A)$ must map $V(L)$ into itself is sufficient
to have the desired irreducibility with respect to the unitary group
generated by $\pi^0(A)$.  We shall leave this as a topic for future
work.

We now look at the case $\kp=1$, that is, $H_{1}=\half(p_x^2+
e^{4x}-p_y^2)$.  Points on $C$ have to satisfy $|p_x|<|p_y|$.  The
flow generated by $H_1$ is, restricted to $C$, $$ (x,y,p_x,p_y)\raw
(x(t),y-p_yt,p_y\tanh[2p_y(t_0-t)],p_y), $$ where $x(t)$ is determined
by the condition $H_1=0$, and $t_0=(2p_y)^{-1}{\rm arctanh}(p_x/p_y)$.
This motion is complete (that is, defined for all $t$), so that
$H_1=J$ is the moment map of an $\R$-action. There are no
singularities, and $S^0$ is duly a manifold, namely $T^*\R$ with the
zero section removed, equipped with the canonical symplectic
structure. In fact, the functions $f_1=p_y$ and $f_2=y-\half {\rm
arctanh}\, (p_x/p_y)$ Poisson-commute with the constraint on $C$, and
project to globally defined canonical co-ordinates on $S^0$.

As before, $\H=L^2(\R^2)$, on which the constraint is quantized by the
closure of the operator $$ \hat{H}_1=\half(-(\partial/\partial x)^2
+\exp(4x)+ (\partial/\partial y)^2),$$ initially defined and
essentially self-adjoint on $C_c^{\infty}(\R^2)$. (The essential
self-adjointness immediately follows from the positivity of
$-(\partial/\partial x)^2 +\exp(4x)$, cf.\ \ci{RS2}, or may be
inferred from the explicit form of the (generalized) eigenfunctions
and the Weyl-Titchmarsh theory, cf.\ \ci{Wei}.)

The operator $\hat{H}_1$ is diagonalized by the unitary transformation
$W_1:\H\raw L^2(\R^+)\ot L^2(\R,dp/2\pi)$ given by $$
(W_1\ps)(\sg,p)=(\ps,f_1(\sg,p;\cdot));\:\:\:
(W_1^{-1}\check{\ps})(x,y)=(\check{\ps},\ovl{f_1(\cdot;x,y)}), $$ with
$$ f_1(\sg,p;x,y)=
\pi^{-1}e^{-ipy}\sqrt{2\sinh(\pi\sqrt{\sg})}K_{i\sqrt{\sg}}(\half
e^{2x}).$$ This is closely related to the (Kontorovich-) Lebedev
transformation, cf.\ \ci{Pic}. The expressions as given are defined on
$\ps\in\H$ and $\check{\ps}\in W\H$ in a suitable dense subset (e.g.,
functions with compact support), and then extended by continuity.  The
point is that $(W_1H_1W_1^{-1}\check{\ps})(\sg,p)=(2\sg-\half
p^2)\check{\ps}(\sg,p)$.  For suitable $L\subset \H$ (as explained
above), the expression $ (\ps,\phv)_0=\int_{\R}dt\,
(e^{-itH_1}\ps,\phv)$ is well-defined, and equal to the inner product
$(V_1\ps,V_1\ps)$ in $\H^0=L^2(\R,dp/2\pi)$, with $V_1$ given by $
(V_1\ps)(p)=(\ps,\tilde{f}_1(p;\cdot)$, where $$
\tilde{f}_1(p;x,y)=
\sqrt{\pi}f_1(p^2/4,p;x,y)=e^{-ipy}\sqrt{\pi^{-1}2\sinh(\pi |p|/2)}
K_{i|p|/2}(\half e^{2x}).$$

For each value of $p$, $\tilde{f}_1(p;\cdot)$ is a solution to the
quantum constraint, i.e., $$\hat{H}_1\tilde{f}_1(p;\cdot)=0.$$ There
is, however, another linearly independent solution as well, which has
been excluded by our method, in the sense that it plays no r\^{o}le in
the construction of the induced space $\H^0$.  Note the difference
with the $\kp=0$ case. There, for each $p$ the quantum constraint had
2 linearly independent solutions as well (viz.\ $f_{\pm}$), which {\em
both} `contributed' to $\H^0$. The discrepancy between the two cases
is a consequence of the fact that the spectrum of $-
(\partial/\partial x)^2$ is $\R^+$ with multiplicity 2, whereas that
of $ -(\partial/\partial x)^2 +\exp(4x)$ is $\R^+$ with multiplicity
1.

Finally, we consider $\kp=-1$, i.e., $H_{-1}=\half(p_x^2-
e^{4x}-p_y^2)$. The peculiar feature of this case is that the flow
generated by $H_{-1}$ is incomplete. Restricting ourselves to the
constraint hypersurface $H_{-1}=0$, we have the following situation.
One has the condition $|p_x|>|p_y|$. For $p_y\neq 0$, the flow is $$
(x,y,p_x,p_y)\raw (x(t),y-p_yt,p_y/\tanh[2p_y(t_0-t)],p_y).  $$ For
$p_x>|p_y|$, this motion is defined for $t<t_0$, and describes how the
$x$-co-ordinate moves from $-\infty$ at $t=-\infty$ to $\infty$ at
$t=t_0$; $p_x$ moves from $|p_y|$ at $t=-\infty$ to $\infty$ at
$t=t_0$. For $p_x< -|p_y|$, the motion is defined for $t>t_0$, and
takes place in the opposite direction. For $p_y=0$, one has $
(x,y,p_x,0)\raw (x(t), y, p_x/(1-\half p_x t), 0)$. For $p_x>0$ the
motion is defined for $-\infty<t<2/p_x$, and for $p_x<0$ it is defined
for $2/p_x<t<\infty$.  Thus the constraint fails to generate an action
of $\R$. Nonetheless, symplectic reduction from this constraint is
well-defined (one follows the `geometric' reduction procedure of
\ci{AGJ}).  The reduced space is a manifold with two components,
symplectomorphic to $T^*\R \cup T^*\R$. To show this, note that
$f_1=p_y$ and $f_2=y-\half {\rm arctanh}\, (p_y/p_x)$ project to
globally defined canonical co-ordinates on the quotient space, the
components of which are projections of the regions $p_x>0$ and
$p_x<0$, respectively.

To quantize, we consider the operator $$
\hat{H}_{-1}=\half(-(\partial/\partial x)^2 -\exp(4x)+
(\partial/\partial y)^2),$$ initially defined on
$D=C_c^{\infty}(\R^2)\subset \H=L^2(\R^2)$.  As can be expected on the
basis of the fact that the classical motion is incomplete
\ci{RS2,Hep}, this operator is not essentially self-adjoint.  The
deficiency indices are $(1,1)$, and the self-adjoint extensions are
characterized by boundary conditions at $+\infty$. For each
$\al\in[0,2\pi)$, the extension $ \hat{H}_{-1}^{\al}$ is the closure
of $ \hat{H}_{-1}$ defined on $D$ with a function added whose
asymptotic behaviour as $x\raw \infty$ is $\sim
z^{-1/2}[\exp(iz)+\exp(-i(z-\al))]$, with $z=\half\exp(2x)$.

It is remarkable that, in general, all self-adjoint extensions of a
given incomplete Hamiltonian have the same classical limit \ci{Hep}.
In the context of constrained quantization, however, we must conclude
that each choice of a self-adjoint quantization of the constraint
leads to a different quantum theory of the reduced system.

For simplicity, we choose $\al=0$. Diagonalization of $\hat{H}_{-1}^0$
(cf.\ \ci{Pic}) is accomplished as in the $\kp=1$ case, the only
difference being that $W_{-1}:\H\raw L^2(\R^+)\ot L^2(\R,dp/2\pi)$ is
now defined by $$ f_{-1}(\sg,p;x,y)= \half e^{-ipy}\sqrt{2{\rm
cosech}\, (\pi\sqrt{\sg})}(J_{i\sqrt{\sg}} +J_{-i\sqrt{\sg}}) (\half
e^{2x}).$$ Consequently, in the definition of $V_{-1}$ one has the
transformation function $$
\tilde{f}_{-1}(p;x,y)= \half e^{-ipy}\sqrt{2\pi 
{\rm cosech}\, (\pi |p|/2)} (J_{i|p|/2}+J_{-i|p|/2})(\half e^{2x}).$$
For each $p$, this is a solution to the quantum constraint.  Since the
spectrum of $ -(\partial/\partial x)^2 -\exp(4x)$ is $\R^+$ without
degeneracy, the same comments as in the $\kp=1$ case apply here.

Here and in all analogous examples, it is possible (though by no means
necessary) to choose $L$ as a Hilbert space $\H_+\subset \H$, in which
case the solutions $f_{(\pm)}(p)$ to the constraints lie in the
continuous dual $\H_-$; the inclusions $\H_+\subset \H\subset \H_-$
are continuous. The `inner product' of $\ps$ and $f(p;\cdot)$ then
stands for the pairing of $\H_+$ and $\H_-$. See \ci{Ber,PSW,PS} for
this approach to self-adjoint differential operators.  Also, it is
possible to choose $\H_+$ as a nuclear space which is continuously
injected in $\H$; this leads to a similar (Gel'fand) triplet
structure. It does not seem to be a good idea to us, however, to
define the physical Hilbert space as the subspace of $\H_-$ consisting
of {\em all} generalized solutions of the constraints \ci{Haj}.  This
space is far too big; firstly it is not possible to make it into a
separable Hilbert space in any reasonable way, and secondly, in the
context of our examples, it contains both linearly independent
solutions of the constraints (for fixed $ p$). While for $\kappa=0$
this happens to be reasonable, for $\kp=\pm 1$ we see that essential
information (on the specific way the constraints are defined as
self-adjoint operators on $\H$) is thrown away. This particularly
affects the so-called `wavefunction of the universe', cf.\ \ci{LanW}.
\section{Mackey's quantization}
To close, we briefly explain how one can understand and generalize
Mackey's quantization method
\ci{Mac} from the point of view of symplectic reduction and Rieffel
induction 
\ci{Lan0,Lan1}.
Mackey studied particle motion on homogeneous spaces $Q=G'/G$ (where
$G$ is a closed subgroup of $G'$), and found that one may define a
family of inequivalent quantizations, one for each element of
$\hat{G}$ (the unitary dual of $G$). These are usually construed as
all being quantizations of the single phase space $S=T^*Q$. This is a
misinterpretation, however. In fact, each is the quantization of a
different phase space.

There is no reason to specialize the discussion to homogeneous spaces;
we will consider a general principal bundle ${\sf P}$ with gauge group
$G$ and base space $Q={\sf P}/G$. The right-action of $G$ on ${\sf P}$
lifts (pulls back) to a strongly Hamiltonian right-action on
$S=T^*{\sf P}$ with moment map $J$ (cf.\ \ci{AM} for an explicit
formula). Given a co-adjoint orbit $\O$, one obtains the
Marsden-Weinstein quotient $\S^{\O}$. As first discussed in physics by
Wong, and mathematically by Sternberg and others, this is the phase
space of a charged particle, moving on the configuration space $Q$,
which couples to a gauge field with group $G$ \ci{Wein}.
  
Now we quantize using our method, for simplicity assuming that $G$ is
compact (which is the case relevant to physics). We equip ${\sf P}$
with a $G$-invariant faithful measure and form $L=\H=L^2({\sf P})$.
This Hilbert space carries a unitary representation $U(G)$ defined by
$(U(g)\ps)(p)=\ps(pg)$.  Also, we assume we have a Hilbert space
$\H_{\rh}$ and a unitary representation $U_{\rh}$ which quantize $\O$
(cf.\ section 2); evidently, this is possible only of the orbit is
indeed `quantizable'. Constructing the modified inner product $(\,
,\,)_0$ on $L\ot\H_{\rh}$ as explained in section 2 (or using the
Dirac approach), one easily finds that the induced space $\H^{\rh}$
may be identified with the (closed) subspace of $L^2({\sf
P})\ot\H_{\rh}$ consisting of equivariant functions, i.e., satisfying
$\ps(pg)=U_{\rh}(g^{-1})\ps(p)$ for a.e.\ $p\in{\sf P}$ and all $g\in
G$. Indeed, this subspace is precisely $P^{\rh}_{\rm id}(L^2({\sf
P})\ot\H_{\rh})$.  (For non-compact $G$, the induced space is no
longer a subspace of $L^2({\sf P})\ot\H_{\rh}$, but is still correctly
identified by Rieffel induction.)  Geometrically, the induced space
$\H^{\rh}$ is the $L^2$-closure of the space of smooth compactly
supported cross-sections $\Gamma_c^{\infty}(E^{\rh})$ of the Hermitian
vector bundle $E^{\rh}$ (over $Q$) associated to ${\sf P}$ by the
representation $U_{\rh}(G)$.

To sum up, the quantization of the Marsden-Weinstein quotient
$S^{\O}$, where $\O$ is a quantizable orbit, is the Hilbert space
$\H^{\rh}$ induced from $\H_{\rh}$ (the quantization of $\O$).  If
${\sf P}$ is a group $G'$ and $G\subset G'$ then this construction
reduces to Mackey's theory of induced representations \ci{Mac}. This
makes it particularly clear how his inequivalent quantizations of the
configuration space $Q$ are to be interpreted.

 A further illustration of the parallel between symplectic quotients
and $C^*$-modules (see section 3) is obtained by noting that the space
${\sf L}=\Gamma_0(E^{\rh})$ of continuous cross-sections of $E^{\rh}$
vanishing at infinity can be made into a $C^*$-module \ci{Con}. The
$C^*$-algebra $\B$ is $C_0(Q)$, which acts on ${\sf L}$ on the right
by $(\pi_R(f)\ps)(q)=f(q)\ps(q)$. The generalized quantum moment map
$\la\cdot,\cdot\ra_{C_0(Q)}:\overline{{\sf L}}\ot {\sf L}\raw C_0(Q)$
is given by $\la\ps,\phv\ra_{C_0(Q)}: q\raw (\phv(q),\ps(q))_{\rh}$.
For any pure state $q\in Q$ on $C_0(Q)$, the Rieffel-induced space
$\H^q$ constructed from these data is just $\H_{\rh}$. (We could have
started from the pre-$C^*$-module ${\sf L}' =
\Gamma_c^{\infty}(E^{\rh})$, with the same conclusion; we used ${\sf
L}=\Gamma_0(E^{\rh})$ because it is complete.)

If our analogy between symplectic reduction and Rieffel induction is
correct, there should be a corresponding construction at the classical
level.  Indeed, we take as ingredients (cf.\ section 1) $S =S^{\O}$
(which is a locally trivial bundle over $Q$ with projection $pr$),
$P=Q$ with zero Poisson structure, and $J:S^{\O}\raw Q$ given by
$J=pr$.  Any point $ q$ (seen as a zero-dimensional symplectic
manifold) in $Q$ defines a symplectic realization $\rh_q:\{q\}\raw Q$
by the inclusion map. The fiber product $C=S^{\O}*_Q \{q\}$ is just
the fiber $pr^{-1}(q)$ above $q$. The null foliation ${\cal F}_0$ is
generated by functions of the type $pr^*f$, $f\in\cin(Q)$.  Locally,
$S^{\O}$ looks like $T^*Q\times \O$ (though not as a symplectic
space), and in this local picture the leaves of ${\cal F}_0$ are the
fibers of $T^*Q$.  Hence the reduced space $S^{\rh_q}$ equals $\O$,
including its symplectic structure.

Turning the argument around, we have written the orbit $\O$ as a
generalized Marsden-Weinstein quotient $S^{\rh_q}$ constructed from
the unreduced space $S^{\O}$.  Assuming that $\O$ is quantized by a
Hilbert space $\H_{\rh}$, we found (from previous analysis) that
$S^{\O}$ is quantized by the linear space ${\sf L}=\Gamma_0(E^{\rh})$.
Applying Rieffel induction, we then showed that the quantization of
the reduced space $S^{\rh_q}$ using our method yields $\H_{\rh}$
again.  Please note that this argument is non-circular.

 \end{document}